\documentclass[12pt]{amsart}

\usepackage{amsthm, amssymb, amsfonts}

\usepackage{mathtext}
\usepackage[cp1251]{inputenc}

\usepackage{url}

\newcommand\R{\mathbb R}

\newcommand{\Cn}{{\mathbb C}^n}

\begin{document}

\title{\bf {Conformality in the sense of Gromov and a generalized Liouville theorem}}

\author{V.A.~Zorich}

\date{April 24, 2021}



\begin{abstract}
M.Gromov extended the concepts of conformal and quasiconformal mapping
to the mappings
acting between the manifolds of different dimensions. For instance, any entire holomorphic function $ f: \Cn \to {\mathbb C}$ defines a mapping conformal in the sense of Gromov.
In this connection
Gromov addressed a natural question: which facts of the classical theory apply to these mappings? In particular is it true that

{\em If the mapping $ F: \R^{n + 1} \to \R^{n}$ is conformal and bounded, then it is a constant mapping, provided that $ n \geq 2 $}~?

We present arguments confirming the validity of such a Liouville-type theorem.
\end{abstract}



\maketitle
\markboth{Vladimir Zorich}{Generalized Liouville theorem}

{\bf 1.~Introduction.}~
\vskip2mm
M.~Gromov proposed the following definition extending the notion of {\em conformal} (quasiconformal)
mapping between domains of Riemannian manifolds of equal dimension.
A mapping between two metric spaces, for instance,
a mapping $ F: \R^{m} \to \R^{n}$ ~ ( ${m} \geq {n}$) is {\em conformal}
if it transforms any infinitely small ball of the domain into an infinitely small ball in the image [1].
A mapping is {\em quasiconformal} if it transforms any infinitely small ball of the domain into an infinitely small ellipsoid whose eccentricity is uniformly bounded  above by a
constant.
For example, any entire holomorphic function
$ f: \Cn \to {\mathbb C}$ defines a mapping conformal in the sense of Gromov.

In connection with such an extension of the concepts of conformality and quasiconformality of a mapping, Gromov addressed a natural question:
which facts of the classical theory could apply to these mappings? In particular is it true that

{\em If the mapping $ F: \R^{n + 1} \to \R^{n}$ is conformal and bounded, then it is a constant mapping, provided that $ n \geq 2 $?}

Below we give arguments confirming the validity of such a Liouville theorem.

\vskip5mm
{\bf 2.~Preliminary observations and discussions.}~
\vskip2mm
To clarify the statement itself and some of its aspects, we start with
some simple observations
useful for further considerations.
\vskip2mm

$ \bullet $ Let us first consider the simplest case, a mapping $ F: \R^{1} \to \R^{1} $.
Differentiability of such a mapping
is, clearly, sufficient for conformality in the sense of Gromov.
However, the Liouville theorem is not applicable
already
to the function
$ x \mapsto y = \arctan x $.

The orthogonal projection $ F: \R^{n} \to \R^{k} $ of the whole space onto its subspace, which includes
the case of a line,
$\R^{k}=\R^{1} $, is a Gromov conformal mapping. The composition of Gromov conformal or quasiconformal mappings, of course, also belongs to this class of mappings.

Hence, for
Gromov conformal or quasiconformal mappings
$ F: \R^{n} \to \R^{1} $ the Liouville theorem does not hold.
\vskip2mm

$ \bullet $
However, the Liouville theorem is applicable to
those functions $ F: \R^{1} \to \R^{1} $ which are not only differentiable, but which also admit a representation by a power series converging on the whole real axis
(since the
Liouville theorem is valid for any entire holomorphic function).
\vskip2mm

$ \bullet $ Finally, if we pass from conformity to quasiconformality,
then it is known that the condition equivalent
to quasiconformality is the so-called \textit{$ h $-condition}.
A mapping $ F: \R^{1} \to \R^{1} $ satisfies the
$h$-condition if for any three points $ a < b
< c $ such that $ |a - b| = |b - c| $ there is a uniform estimate of
the relative distortion of distances:
$$ h^{- 1} \leq \frac{|F (a) - F (b)|}{|F (b) - F (c)|} \leq h \,.
$$
Such mappings arise, for instance, when a
quasiconformal mapping of the half-plane onto itself is extended to
the boundary. Moreover, every mapping satisfying the $ h $-condition
can be extended to a quasiconformal automorphism of the half-plane (as well as
the whole plane). For quasiconformal mappings $ F:
\R^{1} \to \R^{1} $ of that type, the Liouville theorem
holds, similarly to its validity for analytic functions.

Note that the mapping $ F: \R^{1} \to \R^{1} $  satisfying the $ h $-condition
might itself be nondifferentiable at any point but, nevertheless, be quasiconformal
in the classical sense and, thus, in the sense of Gromov as well.
\vskip2mm

$ \bullet $ Now consider the following example of a mapping $ F :
\R^{3} \to \R^{3} $. First, project the entire space $ \R^{3} $ onto
the plane $ \R^{2} $. Then locally isometrically roll the plane into
a cylinder. Then locally conformally (quasiconformally) fold the
cylinder into a torus lying in $ \R^{3} $.
We have constructed a conformal (respectively, quasiconformal)
nonconstant map $ F : \R^{3} \to \R^{3} $
with a bounded image, which shows that the
Liouville theorem does not apply. Note, that
this mapping is not conformal in the sense of Gromov, and
neither is any mapping $ F: \R^{m} \to \R^{n} $ for $ {m} <{n} $.

The above considerations explain why Gromov imposed
conditions $m\ge n\ge 2$ in his question on validity of the
Liouville theorem for
mappings
$ F: \R^{m} \to \R^{n}$ conformal in his sense.

\vskip5mm
{\bf 3.~Preliminary remarks on the proof of the theorem.}~
\vskip2mm
For clarity of interpretation of the arguments carried out below,
we consider the case of a
mapping $ F: \R^{3} \to \R^{2} $. Consider both conformal and
quasiconformal mappings on the same footing.

Let us start again with simple but important
observations, which now deal
not with the formulation but with the proof of the theorem.
\vskip2mm
{\em Some concepts and facts.}
\vskip2mm

$ \bullet $ First of all note that if some mapping
$ F: \R^{3} \to D \subset \R^{2} $ is  generic,
i.e. at each point of the domain the mapping's rank is maximal possible (here it is 2),
then the preimage of any point in $ F (\R^{3}) \subset D $
has to be 
a curve $ \gamma $.

If such a mapping
$ F: \R^{3} \to D \subset \R^{2} $ is
conformal (quasiconformal) in the sense of Gromov, then
its tangent mapping $ F' $,
being restricted to the hyperplane orthogonal to $ \gamma $,
is conformal (quasiconformal) in the usual (classical) sense.

If the distribution of such hyperplanes is integrable,
then, of course, the restriction $ F|_S $ of the mapping $ F: \R^{3} \to D \subset \R^{2} $ to the integral surface $ S $ of the distribution will be conformal (quasiconformal) in the classical sense of conformal or quasiconformal
maps of Riemannian manifolds of
equal dimensions.

We note however that the distribution of hyperplanes induced  by a  mapping
$ F: \R^{3} \to \R^{2} $, even if it is conformal in the sense of Gromov,
does not have to be integrable. This will be demonstrated below in section 5.
\vskip2mm
$ \bullet $ Let us now recall the concept of conformal type of an open Riemannian manifold.
Conformal type is invariant under quasiconformal mappings.

It will be more convenient for us here to define conformal type of a manifold in terms of the conformal modulus of families of curves on the manifold.

Assume that on a Riemannian manifold (on a surface, in our case) there is some family
of curves. Then a nonnegative function on the manifold is said to be
admissible for that family of curves
if the length of each curve of the family weighted by the density defined by that function is at least one. Now integrate
the square of such an admissible function over the entire surface to obtain a weighted surface area. The infimum of the so-weighted surface area over all functions admissible for a given family of curves is called the {\em conformal modulus of the family}.

In the general case of a Riemannian manifold of
dimension $ n \ge 2 $, the weighted volume is obtained by integrating the $n$th power (rather than the square) of an admissible function.

Conformal modulus of the family of curves is invariant under conformal (in the classical sense) mappings of the  Riemannian manifold.
Under  a quasiconfor\-mal mapping it changes by a factor
depending on the
coefficient of quasicon\-for\-ma\-li\-ty of the mapping.

For conformal mappings that are not homeomorphisms, the conformal modulus may change, but only by getting smaller.
Similarly, for nonhomeomorphic
quasiconformal mappings one can also guarantee that if the modulus of the preimage is equal to zero, then the modulus of the image cannot be positive, and hence must be equal to zero.

In principle, this information would be sufficient for our goal, but it is still useful to say a few words about the conformal type of an open connected Riemannian manifold.
On such a manifold take a compact set that does not degenerate into a point, for example, take a geodesic ball
--- a neighbourhood of some point. Consider a family of all curves starting
at the compact set and going to infinity in the sense that they go beyond any compact part of the manifold.
A manifold is said to be
{\em conformally parabolic} or, respectively, {\em conformally hyperbolic}
if
the modulus of such a family is zero or, respectively, positive.

For example, the Euclidean space $ \R^{n} $ is conformally parabolic, while the Lobachevsky space
$ {\mathbb H}^{n} $ is conformally hyperbolic.
Just for this reason,
they are not conformally equivalent and there is no holomorphic mapping of a plane into a circle (disk) different from a constant.
(One should not associate the conformal type of a surface or a Riemannian manifold with its curvature. The ordinary curvature of a surface is not a conformal invariant, as the stereographic projection of the sphere onto  the plane or the Poincar\'e model of the Lobachevsky plane in a disk show.)

\vskip2mm

$ \bullet $ What may happens to the conformal type of a surface when the mapping is conformal in the sense of Gromov?

Imagine a surface $ H $ in $ \R^{3} $ which is a graph of a function defined
on the whole plane $ \R^{2} $. The plane itself is a conformally parabolic surface,
but the graph may happen to be a surface of conformally hyperbolic type (although the corresponding construction is nontrivial,
see, for example, [2]). 

By the orthogonal projection $ F: \R^{3} \to \R^{2} $, which is obviously conformal in the sense of Gromov,
the surface $ H $ of hyperbolic type is transformed into the plane $ \R^{2} $.

Thus, conformal mappings in the sense of Gromov, generally speaking, do not preserve the usual conformal types of
objects.
In the example considered above the point is, of course, that the restriction of the projection mapping $ F: \R^ {3} \to \R^{2} $ onto the surface $ H $ is neither  conformal, nor
 quasiconformal in the classical sense of a mapping of the surface $ H $ to the plane $ \R^{2} $.

But if some surface $ S $ turns out to be an integral surface of the distribution induced by the mapping $ F: \R^{3} \to \R^{2} $, conformal or quasiconformal in the sense of Gromov, then
the restriction $ F|_S $ of the mapping $ F $ to the surface $ S $ is, respectively, conformal or quasiconformal in the classical sense.
The conformal type of such a surface is preserved even under a Gromov generalized quasiconformal mapping.

One can find more useful information about the conformal types of Riemannian manifolds including geometric criteria
of conformal type of a manifold, for example, in [3], [4], [5].

\vskip5mm
{\bf 4.~Description of the main constructions.}~
\vskip2mm
Now we assume that we are dealing with a mapping $ F: \R^{3} \to D \subset \R^{2} $
nowhere degenerate and sufficiently regular (smooth).

\vskip2mm

$ \bullet $ For each compact ball
centered at the origin of the space $ \R^3 $
there exists a constant $ \delta = \delta (r) $, depending only on the radius of the ball, such that
at all points of the ball whose mutual distance
does not exceed
$ \delta = \delta (r) $, normals to the planes of distribution
generated by the mapping $ F: \R^{3} \to D \subset \R^{2} $
have angle deviation from each other bounded above by some fixed small value $ \varepsilon> 0 $.

\vskip2mm

$ \bullet $ We
would like to lift paths from the image $ F (\R^{3}) \subset D \subset \R^{2} $ to the covering space
$ \R^{3} $. For this we first say how the vector fields
lift to the covering space.
Let
$ \xi \in F (\R^{3}) $ and $ \tilde \xi $ be a point in the space $ \R^{3} $ such that
$ F (\tilde \xi) = \xi $.
The plane of the distribution, generated in the space $ \R^{3}$ by the nowhere degenerate map $ F $ associated with the point $ \tilde \xi $, under the action of the nondegenerate tangent map $ F '(\tilde \xi) $
linearly
transforms to the two-dimensional space $ T_{\xi} \R^{2} $ tangent to the image $ F (\R^{3}) \subset \R^{2} $ at the point $ \xi $. Hence, to each vector $ v \in T_{\xi} \R^{2} $ we associate a well-defined vector
$ \tilde v $ attached to the point $ \tilde \xi $.

Thus, vector fields distributed in the domain $ F (\R^{3}) \subset \R^{2} $ admit lifting to the space $\R^{3}$.

It is now clear that any smooth path $ \gamma \subset F (\R^{3}) $ emanating from the point $ \xi $
admits a unique
covering path $ \tilde \gamma \subset \R^{3} $ emanating from the point $ \tilde \xi $
and everywhere tangent to the planes of the distribution
generated in the space by the mapping $ F $.

\vskip2mm

$ \bullet $
Let us make one more observation, obvious but very useful,
concerning the lifting of
the path from the image to the preimage.
Take
a point
$ \xi \in F (\R^{3}) \subset D $
and some of its preimages $ \tilde \xi \in \R^{3} $.
In the domain $ D $ consider a
ray $ R_{\xi} $ emanating from the point $ \xi $. In the preimage, to the path $ R_{\xi} $ one associates
the path $ \tilde R_{\xi} $ starting at the point $ \tilde \xi $ and going
along the planes of the distribution induced by the map $ F: \R^{3} \to D \subset \R^{2} $.

It is important that as the path $ R_{\xi} $ approaches the boundary of the domain $ D $
(more precisely,  the boundary of $ F (\R^{3}) \subset D $), the path $ \tilde R_{\xi} $
has to go to infinity in the space $ \R^{3} $.
We will heavily use this fact below.

\vskip2mm
$ \bullet $
Finally here is one cautionary note related to the possible nonintegrability
of the distribution of hyperplanes generated in the space $ \R^{3} $ by the map
$ F: \R^{3} \to D \subset \R^{2} $.

In the domain $ F (\R^{3}) $ take a rectangle spanned by a horizontal segment $ I $ and a height $ h $.
Let $ \xi_0 $ and $ \xi_1 $ be the start and end points of the segment $ I $.
One can reach the vertex of the rectangle opposite to the vertex $ \xi_0 $
in two ways. One can first go horizontally and then vertically, or one can first go up vertically and then go horizontally.
Taking some preimage $ \tilde \xi_0 $ of the point $ \xi_0 $, each of these two paths, starting at the point $ \tilde \xi_0 $, can be lifted
into the space $\R^{3}$. But
the ends of the lifted paths
may not coincide now. This is the effect of the nonintegrability of the distribution.
This remark
partly clarifies the meaning of the construction below.

\vskip2mm
$ \bullet $ Let $ \tilde I \subset \R^{3} $ be a path covering the horizontal segment $ I \subset \R^{2} $ selected above and starting at point $ \tilde \xi_0 $.
At
each point $ \xi $ of the segment $ I $ take
a ray $ R_{\xi} $ going vertically upward.
Its covering path $ \tilde R_{\xi} $ starting at $ \tilde I $ is now uniquely determined.
(Of course, we can lift only that part of the ray $ R_{\xi} $ that is contained in $ F (\R^{3}) $.)

In particular, consider
the vertical ray $ R_{\xi_0} $ and let $ \tau $ be the height of the point of the ray
$ R_{\xi_0} $ over (i.e. above) the point $ \xi_0 $. We will identify the number $ \tau $ and the corresponding point of the ray $ R_{\xi_0} $.

\vskip2mm
$ \bullet $ Consider the initial segments of the rays $ R_{\xi} $ of small height $ h $.
Their union (for all $ \xi \in I $)
sweep out the rectangle $ S_h^R $ with base $ I $ and the height $ h $.
If for each value  $ \tau \in [0, h] $ we take the segment $ I_ \tau $ obtained from the segment
$ I =: I_0 $ by the vertical shift by $ \tau $, then the segments $ I_ \tau $ will sweep the same rectangle
$ S_h^R $, which we will denote this time by $ S_h^I $.

\vskip2mm
$ \bullet $
The curves
$ \tilde R_{\xi} $ that are lifting of pieces of rays $ R_{\xi} $ lying in the rectangle $ S_h^R $
sweep out some surface $ {\tilde S}_h^R \subset \R^3 $.

The curves
$ {\tilde I}_\tau $, the lifting of the segments $ I_\tau $
for $ \tau \in [0, h] $, sweep out some
surface $ {\tilde S}_h^I \subset \R^3 $.

Despite the fact that $ S_h^R = S_h^I $, the surfaces $ {\tilde S}_h^R $ and $ {\tilde S}_h^I $, as we have already noted, may not coincide if the distribution of hyperplanes generated in the space $ \R^{3} $ by the map
$ F: \R^{3} \to D \subset \R^{2} $ is nonintegrable.

\vskip2mm
$ \bullet $ Using the
smoothness
of the distribution of hyperplanes, we now choose a sufficiently small positive number $ h_1 $
so that the restriction of the mapping $ F $ to the surface
$ {\tilde S}_{h_1}^R $
be
quasiconformal with the coefficient of quasiconformality, for example, not greater than twice of the quasiconformality coefficient of the mapping $ F $ itself.

\vskip2mm
$ \bullet $ Thus, starting from the segment $ I = I_0 $ and its cover
$ \tilde I $, lifting pieces of vertical rays $ R_ \xi $ of small height $ h_1 $, we covered the surface
$ {\tilde S}_{ h_1}^R $. The image of this surface
is a rectangle with horizontal sides $ I = I_0 $ and $ I_ {h_1} $.

It is important that the restriction of the mapping $ F $ to the constructed surface $ {\tilde S}_{h_1}^R $ is of
controlled coefficient of quasiconformality, which
is at most twice of the quasiconformality coefficient of the mapping $ F $ itself.

\vskip2mm
$ \bullet $ We have described the first step of the construction. Denote $ {\tilde S}_{h_1}^R $ by a shorter symbol $ {\tilde S} _1^R $.

Now, taking instead of $ I = I_0 $ and $ \tilde I = \tilde I_0 $
segment $ I_ {h_1} $ and its lifting $ {\tilde I}_{h_1} $,
repeating the same procedure, and shifting in height
from $ h_1 $ to $ h_2 $, we
construct the surface $ {\tilde S} _2^R $.

Note
that the surface $ {\tilde S} _2^R $, generally speaking, is not a continuous conti\-nu\-ation of the surface $ {\tilde S} _1^R $.
Together they rather resemble the steps of a staircase connected by a railing which in this case is
the lift $ \tilde R_{\xi_0} $ of the ray $ R_{\xi_0} $.

\vskip2mm
$ \bullet $ What could prevent further expansion of the surface
$ \tilde S = \bigcup \tilde S_i^R $ constructed by the  procedure described above?

\vskip2mm
The only obstruction to continuation of this construction may
occur in the following way.
For some limiting value $ h_ \infty $ in the segment $ I_ {h_ \infty} $ there appears a singular point of
such a kind that while approaching this point in the image
the  corresponding inverse images tend to
infinity in the space $ \R^{3} $.

To avoid  unnecessary technicalities, let us first assume that the entire segment
$ I_ {h_ \infty} $ consists of such singular points that prevent expanding the
surface $ \tilde S $.

\vskip2mm
$ \bullet $
The covers of the rays $ R_{\xi} $, starting in the segment $ I $ and going to the points of the segment $ I_ {h_ \infty} $, are the step curves $ \tilde R_{\xi} $ lying on the staircase surface
$ \tilde S = \bigcup \tilde S_i^R $
and going
to infinity in the space $ \R^{3} $.

The conformal modulus of the
family $ R = \{R_{\xi} \} $ of rays $ \{R_{\xi} \} $ is positive, since the domain $ D $ is supposed to be
bounded.
(Actually this conformal modulus is equal to the conformal modulus of
the rectangle with the parallel sides $ I $ and $ I_ {h_ \infty} $.)

If we show that the family $ \tilde R $, the preimage
of the family $ R $
in the surface $ \tilde S $, is of zero conformal modulus ($M_2 (\tilde R) = 0$),
then we come to
a contradiction, since
the restriction $ F|_{\tilde S} $ of the mapping $ F: \R^{3} \to D \subset \R^{2} $ to the surface $ \tilde S $,
by construction,
is quasiconformal in the classical sense.

\vskip2mm

$ \bullet $
To verify the equality  $M_2 (\tilde R) = 0$,
let us first show  that if the mapping $ F: \R^{3} \to D \subset \R^{2} $ is quasiconformal in the sense of Gromov, then,
roughly speaking,
the area of the part of the surface $ \tilde S $, lying in the ball of radius $ r $
centered at the origin of the space $ \R^{3} $, grows no faster
than $ O (r^2) $ as $ r \to \infty $.

Let us slightly clarify that statement.

Take the inverse image $ \tilde I $ of the segment $ I $ from which the surface $ \tilde S $ began to grow.
In the inverse image $ \tilde R_\xi $ of each ray $ R_\xi $ we take only the part of
length
$ r $
counted from
$ \tilde I $.
Since the inverse images $ \tilde R_ \xi $ of the rays $ R_\xi $ go to infinity, then for large values
of $ r $ the set $ \tilde I $ in the described situation will be equivalent to a reference point.

For large values  $ r $
the area of the `disk' ~of radius $ r $ centered at  $ \tilde I $ grows in the surface $ \tilde S  $ as $ r \to \infty $.

The assertion is that if the mapping $ F: \R^{3} \to D \subset \R^{2} $ is quasiconformal in the sense of Gromov, the  area mentioned grows no faster than $ O (r^2) $.

Since the increment of the area of
a disk can be written in the form $ L (r) dr $, where $ L (r) $ is the length of a `circle'
of radius $ r $, it suffices to check that the length of the circle of radius
$ r \to \infty $ grows
no faster than $ O (r) $.

\vskip2mm

$ \bullet $
Otherwise  the following situation is inevitable, since the surface lies in the Euclidean space:
steep waves should appear on this circle (`steep' in the sense that their height is much greater than their length).
Thus,
the surface $ \tilde S $ creates
something like a gorge  in the space $ \R^{3} $. Namely, there arise  two
almost parallel pieces of the surface $ \tilde S $, that are opposite of each other and are at a distance arbitrarily small in comparison with the height of the gorge.
But the surface $ \tilde S $ is almost integral in the sense that its tangent planes are of nearly the
same direction as the planes of
the distribution generated by the mapping $F$ at the corresponding points of the space.
Hence, the lines $ \gamma $,
along which the mapping $ F: \R^{3} \to D \subset \R^{2} $ projects the space into the plane, being
orthogonal to the planes of the distribution, turn out to be almost orthogonal to the surface $ \tilde S $.
The lines $ \gamma $ must not intersect (otherwise they coincide).
But the curves $ \gamma $
should not coincide since the opposite pieces of the surface have different images in the plane.
Since two almost parallel opposite pieces of the surface $ \tilde S $ are very close and the  curves $ \gamma $ passing through them are almost orthogonal to the surface $ \tilde S $,
the curves $ \gamma $ should create two flows in the gorge.
Moving away from the corresponding piece of the surface in almost normal directions, they immediately
change its direction to parallel to the surface in order to exit
from a narrow crevice.
The density of the curves at the exit from this crevice
exceeds the density of their distribution over gorge's walls.
The ratio of the densities grows in accordance with the ratio of the height of the gorge and its width.

Now let us consider a small three-dimensional ball in the place where the density of curves $ \gamma $  at the exit from the crevice is high, and let us look at the image of this ball under for the initial map
$ F: \R^{3} \to D \subset \R^{2} $.

Let us show that the image of the radius of the ball, which is directed like a bridge across the crevice,
will be much larger than the image of the radius directed
along the crevice.

Indeed, consider the projection of those radiuses
 by the curves $ \gamma $ onto the surface of the gorge.
The transverse radius is projected along the vertical direction of the gorge.
Its length increases very significantly, since the curves $ \gamma $ are almost parallel
to the lateral surface of the gorge.

On the other hand, the longitudinal radius is projected into a vector of almost the same longitudinal direction.
Its length almost does not change.

Thus, the eccentricity of the image in the surface $ \tilde S $ of the considered ball
can be made arbitrarily large.
But by the construction of the surface $ \tilde S $ the restriction $ F|_{\tilde S} $ of the mapping $ F: \R^{3} \to D \subset \R^{2} $
onto the surface $ \tilde S $ is quasiconformal in the classical sense.

This contradicts to  the assumption that
the initial mapping $ F $ was quasicon\-for\-mal in the sense of Gromov.

\vskip2mm

$ \bullet $
Now let us return to the family of curves $ \tilde R $ in the surface $ \tilde S $ and show that its conformal modulus is zero.
Indeed, since the curves of the family go to infinity, then
the function $ \rho = \alpha ~ (r \ln r)^{- 1} $ for $ \alpha> 0 $ is
admissible for curves of this family, since $ \int_{\gamma} \rho = \infty $ for any
curve $ \gamma = \tilde R_\xi \in \tilde R $. On the other hand, the integral of the square of this function
over the part of the surface that lies in a neighbourhood of the infinity of the space $ \R^3 $
is finite and tends to zero as $ \alpha \to 0 $. Hence, indeed, $M_2 (\tilde R) = 0$.

(We used here the asymptotics $ L (r) = O (r) $ of
the length  of the geodesic circle
of radius $ r $ for $ r \to \infty $.
It is worth mentioning that, according to the criteria proved in [4],
the surface $ \tilde S $ is a surface of conformally parabolic type.)

\vskip2mm

$ \bullet $
It remains
to remove the simplifying technical assumption made above that the entire segment
$ I_{h_\infty} $ consists of singular points
preventing further extension
of the
surface $ \tilde S $.

If there is a point on the segment $ I_{h_\infty} $
which does not prevent expanding
the surface $ \tilde S $,
then
all points in some   neighbourhood $ U \subset I_{h_\infty} $ are nonsingular in the same sense.

Now  the segment $ U \subset I_{h_\infty} $
can play the role
of the initial segment $ I $.

Of course, one needs
to check
that if we take two segments
$ I' \subset U $ and $ I'' \subset U $ having a
common point, then
 continuations of the surface $ \tilde S $ through $ I'$ and through $ I''$ near the common part
$ I' \bigcap I'' $ of these segments will match. This verification is straightforward.

Iterating the described procedures of continuation through nonsingular points we obtain
 a strip-like region in the image $ F (\R^3) \subset D \subset \R^2 $. One part of its boundary
will be the initial
segment $ I $ and, instead of the boundary segment $ I_{h_\infty} $, now there will be a set of points that are singular in the sense
that one cannot  continue lifting through them the rays $ R_\xi $ that form the strip.

Now one can repeat without any changes the
reasoning related to the comparison of the conformal modulus
of the ray family $ R $ (which is not equal to zero provided that the domain $D$ is bounded), and the
the conformal modulus of the family $ \tilde R $ of their lifts (which is equal to zero).

\vskip2mm
$ \bullet $
So, by assuming that the mapping $ F: \R^3 \to D \subset \R^2 $ is quasiconformal in the sense of Gromov, is nonconstant, and its image is a bounded domain $ D $, we arrive at a contradiction.

Hence, indeed, the following Liouville-type theorem holds.
\vskip2mm

{\bf Theorem}

{\em
There is no nonconstant mapping $ F: \R^3 \to D \subset \R^2 $ quasiconformal in the sense of Gromov whose image is a
 bounded domain $D$.
}

\vskip5mm
{\bf 5.~Some additions and comments.}~
\vskip2mm
{\em On the concept of conformity.}
\vskip2mm
$ \bullet $  Conformal mappings are usually discussed in connection with
the geometric meaning of modulus and argument of the derivative of a holomorphic function $ w = f (z) $.
 Note that a nondegenerate linear map $ L: \R^n \to \R^n $ preserves angles if and only if it is a dilation, up to a space motion. For this reason, speaking about the conformity of a mapping on the plane, it is sufficient to know that it transforms infinitesimal circles into infinitesimal circles (infinitesimal balls into infinitesimal balls in higher dimensions).

In terms of conservation of infinitesimal circles (balls) the concept of conformity is meaningful only for
mappings of spaces endowed with suitable metrics, e.g. a Riemannian one.

This remark also applies to the concept of quasiconformal mapping. Con\-for\-mal and quasiconformal mappings
are usually considered for domains in spaces of equal dimensions. But the described constraint
that an infinitesimally small ball is transformed into an infinitesimally small ball (or to a distorted ball) is also applicable when
we are dealing with mappings of spaces of different dimensions, e.g.
a mapping $ F: \R^3 \to \R^2 $.

This is precisely what one finds in the generalized concept of conformality (respectively, quasiconformality) of a mapping in the sense of Gromov.

\vskip2mm
$ \bullet $
There are  several `Liouville theorems' related to various objects.
One of the classical Liouville theorems asserts that all con\-for\-mal mappings in the spaces
$ \R^n $ of dimension $ n> 2 $  are exhausted by those of the M\"obius group (generated by isometries, dilations, and inversion). In this sense the space domains are conformally rigid. For example, it is impossible to conformally map a sphere onto an ellipsoid, unless the latter is a sphere itself.

But one should distinguish between conformal mappings and conformal changes of metric. For example, the Poincar\'e model
of the Lobachevsky plane (space) is obtained by a conformal
change of the Euclidean metric
in a circle (ball). A conformal change of a metric is
a local change of the scale of measurement of lengths,
and an infinitesimal sphere remains
a sphere with respect of the new metric.

Recall that not only conformal mappings, but even
conformal changes of a metric, preserve the conformal type of noncompact Riemannian manifolds.
(Details can be found in [4], [5]). For instance, it is impossible to obtain the Euclidean plane from the Lobachevsky plane by a conformal change of its Riemannian metric (and vice versa). And this statement is quite different from a Liouville-type theorem that there is no bounded entire holomorphic function other than a constant.

\vskip2mm
{\em On the specifics of mappings conformal in the sense of Gromov.}

\vskip2mm
$ \bullet $
We know that the restriction of the mapping $ F: \R^3 \to \R^2 $ conformal or quasiconformal in the sense of Gromov to a two-dimensional surface lying in the space $ \R^3 $ does not have to be conformal or quasiconformal
in the classical sense. In particular, this is why the conformal type of a surface may change under such mappings.

For instance, take the domain $ \R^3 \setminus \R$ where
$\R$ is a straight line.
Consider the domain $\R^2 \setminus O $, the punctured plane,
and take the orthogonal projection
$ \R^3 \setminus \R \to \R^2 \setminus O $
along the line $\R$.
This mapping is conformal in the sense of Gromov and it transforms
the conformally hyperbolic object
$ \R^3 \setminus \R $ into the object $ \R^2 \setminus O $ of conformally parabolic type.

\vskip2mm
$ \bullet $
Note that the distribution of planes induced by the mapping $ F: \R^{3} \to \R^{2} $ quasiconformal or even
conformal in the sense of Gromov
can indeed be nonintegrable.

Let us recall the Hopf bundle $ S^3 \to S^2 $. Making inversions at the poles of the spheres $ S^3 $ and $ S^2 $
get the mapping $ \R^3 \to \R^2 $ conformal in the sense of Gromov.
The corresponding distribution is nonintegrable. The mapping $ \R^3 \to \R^2 $ is, of course,  unbounded.

\vskip2mm
$ \bullet $
As a slightly simpler example of a nonintegrable distribution of hyperplanes in the space $ \R^3 $ we can consider the distribution of hyperplanes orthogonal to helical lines with a fixed pitch
along the vertical axis $ Z $. As a mapping, consider the natural projection of $ \R^3 \setminus Z $ along these helical lines onto a half-plane
whose vertical edge lies on the $ Z $ axis.

The quasiconformality coefficient of this mapping increases when approaching the $ Z $ axis. Outside the vicinity of the axis the mapping is quasiconfor\-mal, and it becomes more and more conformal
as one moves away from this axis

Away from the $ Z $ axis the vertical planes containing this axis are almost integral surfaces of the described distribution. They are certainly of the conformally parabolic type.

Perform the inversion of the space.
Near the origin, the image of infinity, the image of the distribution is arranged as follows.
Along the $ Z $ axis, near the origin, the planes of the distribution are orthogonal to the $ Z $ axis, but in directions orthogonal to the $ Z $ axis the distribution planes
are parallel to the $ Z $ axis. The origin is a singular point.
There are almost orthogonal planes of distribution which are, nevertheless, arbitrarily close to each other.
But their projections along the helical lines do not intersect. And in the initial space
$ \R^3 \setminus Z $ the corresponding points and planes are
far from each other.

\vskip2mm
{\em Integrability and connectivity.}
\vskip2mm
$ \bullet $
In mechanics, where nonholonomic constraints arise (for example, when describing the movement of a skate),
or in classical thermodynamics (when describing adiabatic transitions between states of a thermodynamic system),
the following mathematical problem arises. Consider a distribution of hyper\-plan\-es in space and
allow motions only along the curves integral for this distribution. Under what conditions on the distribution the
space remains connected in the sense that from any of its point one can go to any other point subject to the indicated constraint that one moves
only along the hyperplanes of the distribution?

The answer is given by the  Caratheodory theorem, which he proved (following Poincar\'e)
when considering  the mathematical aspects of classical thermodynamics. Roughly speaking, the theorem claims
that the absence of the connectivity appears only if the distribution is integrable.

\vskip2mm
$ \bullet $
In this regard let us pay attention to the following phenomenon.

Given a nowhere degenerate vector field  in the Euclidean space $ \R^3 $,
it generates a distribution of hyperplanes orthogonal to vectors of this field.

We know that if a vector field is potential, then the work of the field along any path with a fixed start and a fixed end is the same, it does not depend on the path. The distribution of planes associated with the potential field is obviously integrable. Its integral surfaces are levels of the potential.

Conversely, if the initial vector field is nonpotential (more precisely, when it is  not potential anywhere, even locally), then the corresponding distribution of planes is com\-ple\-tely nonintegrable.

In such a field you can do the following. Given any path $ \gamma $ from $ A $ to $ B $,
in arbitrarily small neighbourhood of that path
there exists another
path $ \gamma_0 $ from $ A $ to $ B $ that is integral for the distribution, and hence
the work of the field along this path is equal to zero.
After the Archimedes lever (and defining the elementary work as the scalar product of the force and the displacement),
we know that the gain in effort is achieved by the loss in the path length.
In our problem, this boils down to different path lengths for various paths
between fixed states  in a potential field.
If the field is not potential, then, as we see now, it is possible to move from one state to another even with zero
work of the field.

It is clear that this is a property of any path going along the planes of the distribution, since
the path is orthogonal to the vectors of the field.{\footnote {As a byproduct,  let us build
a perpetuum mobile.
In the standard Euclidean space $ \R^3 $ with Cartesian coordinates $ (x, y, z) $
take the differential 1-form $ \omega = dz $ and perturb it a little to obtain the form
$ \tilde \omega = dz + \varepsilon ~ xdy $.
The 1-form $ \omega $
defines the integrable distribution of the planes $ \ker dz $ in space, while the distribution of the planes
$ \ker (dz + \varepsilon ~ xdy) $ given by the 1-form $ \tilde \omega $,
on the contrary, is completely nonintegrable.
Consider the fields $ V $ and $ \tilde V $ of vectors orthogonal to the planes of
each of these two distributions, respectively.
In the first case the vector field $ V $ is potential with the potential $ z $. We can treat $ V $ as a local gravitational field and $ z $ as a height. If the unit mass is lifted from the origin
along the vertical $ \gamma = [A, B] $ from the level $ z = 0 $ to the level $ z = 1 $, it will have some potential energy which it can release when falling from level $ z = 1 $ to level $ z = 0 $.
Now take the path $ \tilde \gamma $ that goes from $ A $ to $ B $ along the planes of the distribution induced by the form
$ \tilde \omega $.
Since this path goes in the direction perpendicular to the vectors of the field $ \tilde V $, the work of the field $ \tilde V $ on such a transition from $ A $ to $ B $ will be equal to zero.
Now consider the field $ \tilde V $ alone.
Let a small bead of unit mass  fall from $ B $ to $ A $ along a smooth vertical thread and then return from $ A $ to
$ B $ along $ \tilde \gamma $, then we get a perpetual motion machine. (Let us call it the Petrik engine).
}}

\vskip2mm
{\em Towards the Liouville theorem on a constant.}
\vskip2mm
$ \bullet $
The classical Liouville theorem, present in any course of complex analysis, asserts that if an entire holomorphic function is bounded, then it is constant.
The proof is usually given in terms of the Cauchy integral representation of  coefficients of the power series expansion of
 this function.

But such a theorem, of course, also holds for harmonic functions, and moreover not only in the plane,
but in the space as well.

Let the harmonic function be understood as a steady-state temperature distribution in the domain for a given temperature regime
at its boundary. Then naturally the temperature should be vary very little not only between any  close points
of the domain, but also for any two fixed points  far from the boundary, provided that the temperature at the boundary or the flow of energy through the boundary is uniformly bounded.
Hence by moving the boundary to infinity we obtain the Liouville theorem, it is as simple as that!

\vskip2mm
$ \bullet $
Hence, such a Liouville theorem holds for any class of functions and mappings, provided that there is a suitable connection between the behaviour of the mapping inside the domain and its behaviour on the boundary. For example, it is the case
if the maximum principle holds, which means the openness of the mapping, accompanied by the fact of
weakening (averaging) of the influence of the boundary conditions as we move from the boundary into the domain.
\vskip2mm
The Liouville constant theorem for holomorphic or harmonic functions is, of course, a special case
of this more general phenomenon, just as,
for instance, the fundamental theorem of algebra on roots of a polynomial is a theorem not about a polynomial per se,
but about the index of a mapping  applied to the the mapping $ z^n + o (z^n) $.

Even the classical Liouville theorem on  constancy of a bounded entire function is, in fact, related only to the behaviour of a holomorphic function in a punctured neighbourhood of a singular point (infinity in that case). Indeed, if the function is bounded in such a neighbourhood, then it has a finite limit there. For an entire function this leads to the
violation of the mapping's openness. And if only a punctured neighbourhood of an isolated singular point is considered, then  the Picard theorem already prohibits this point from being an essential singularity of a holomorphic function, provided that
the values of that function avoid some set containing more than two points, not to mention a disk or its exterior.

Note that in our reasoning carried out above and
proving the Liouville theorem for mappings quasiconformal in the sense of Gromov,  we actually studied only
the behaviour of the mapping in a neighbourhood of infinity.

\vskip2mm
{\em Towards one old problem.}
\vskip2mm
In [7] we recalled, in particular, a question that arose shortly after the global homeomorphism theorem for quasiconformal mappings was proved [8]. The theorem claims that

{\em Any locally invertible quasiconformal mapping of the space $ \R^n $ into itself is globally invertible provided the dimension $ n> 2 $.}

The remaining question was if this statement also holds
in the infinite-dimen\-si\-o\-nal case? Is it true that {\em if a nonlinear
operator acting in a Hilbert space is locally invertible, then, under the condition of quasiconformality, it is globally invertible?}
The reasoning given above in the proof of Liouville theorem for mappings quasiconformal in the sense of Gromov can be used to confirm an affirmative answer to that question.

Regarding the theorem proved above, we note in conclusion  that while
we considered only the case of a mapping $ F: \R^{3} \to \R^{2} $, in fact
the theorem holds in the general case for a mapping $ F: \R^{m} \to \R^{n} $ for $ m \geq n \geq 2 $. Indeed,
keeping the condition of the quasiconformality of the mapping in the sense of Gromov, we can go
from the mapping $ F: \R^{m} \to \R^{n} $
 to the case $ F: \R^{m} \to \R^{2} $ by projection  and apply the same construction.

\vskip2mm
{\em Final remark.}
\vskip2mm

In note [6], we have already touched upon the main issue considered here.
There we indicated analytic conditions for the conformality of the mapping $ F: \R^{3} \to \R^{2} $. They correspond to
the classical Cauchy--Riemann conditions. In [6] we also
made some observations that turned out to be useful and were implemented above.

For the information of the reader we also note
that on the very first page of the article [1], in a footnote, Gromov gives the address where a much more complete text of his work in English is presented.
In particular,  the reader can find an extended interpretation of conformity and quasiconformality there, as well as
the formulation of the question on a Liouville-type theorem for such mappings.
Now that text is available at
\small {\url {https://www.ihes.fr/~gromov/wp-content/uploads/2018/08/problems-sept2014-copy.pdf}}.

\vskip7mm

\centerline {\bf REFERENCES}
\vskip2mm
\vskip2mm

[{\bf 1}] M.L.~Gromov, “Colorful categories”, Russian Math. Surveys 70:4 (2015), 591-655.

[{\bf 2}] {R.~Osserman}, “Hyperbolic surfaces of the form $ z = f (x, y)” $. Math. Ann. 144 (1961), 77-79.

[{\bf 3}] {J.~ Milnor}, “On deciding whether a surface is parabolic or hyperbolic”.
Amer. Math. Monthly, 84:1 (1977), 43-46.

[{\bf 4}] V.A.~Zorich, V.M.~Kesel'man, “On the Conformal Type of a Riemannian Manifold,” Funct. Anal. Appl., 30:2 (1996), 106-117.

[{\bf 5}] V.A.~Zorich, V.M.~Kesel'man, “The Isoperimetric Inequality on Manifolds of Conformally Hyperbolic Type”, Funct. Anal. Appl. 35:2 (2001), 90-99.

[{\bf 6}] V.A.~Zorich, “On a question of Gromov concerning the generalized Liouville theorem”,
Russian Math. Surveys, 74:1 (2019), 175-177.

[{\bf 7}] V.A.~Zorich, “Some observations concerning multidimensional quasiconformal mappings”, Sb. Math., 208:3 (2017), 377-398.

[{\bf 8}] V.A.~Zorich, “The global homeomorphism theorem for space quasiconformal mappings, its developement and related open problems”, Lecture Notes in Math., 1508 (1992), 132-148.

\vskip5mm

{\small {\bf V.A.~Zorich}

Lomonosov Moscow State University,

\vskip2mm
E-mail: {\bf vzor@mccme.ru}
}

\end{document}